\documentclass[12pt]{amsart}
\usepackage{amsfonts, amsmath}
\newtheorem{theorem}{Theorem}[section]
\newtheorem{prop}[theorem]{Proposition}

\newtheorem{lemma}[theorem]{Lemma}

\newtheorem{cor}[theorem]{Corollary}

\def\Ma { {\mathcal M} }

\newcommand\R{\mathbb{R}}

\title[ Weighted norm inequalities and sharp spectral multipliers]  
{ Weighted norm inequalities, Gaussian bounds and sharp spectral multipliers }

\medskip

\author{Xuan Thinh Duong, Adam Sikora and Lixin Yan}
\address {Xuan Thinh Duong, Department of Mathematics, Macquarie University, NSW 2109, Australia}
\email{xuan.duong@mq.edu.au}
\address{
Adam Sikora, Department of Mathematics, Macquarie University, NSW 2109, Australia}
\email{
adam.sikora@mq.edu.au}
\address{
Lixin Yan, Department of Mathematics, Sun Yat-sen (Zhongshan) University, Guangzhou, 510275, P.R. China}
\email{
mcsylx@mail.sysu.edu.cn
}

\begin{document}

\begin{abstract}
Let $L$ be a non-negative self adjoint operator acting on $L^2(X)$
where $X$ is a space of homogeneous type. Assume that $L$ generates
a holomorphic semigroup $e^{-tL}$ whose kernels $p_t(x,y)$
have Gaussian upper bounds but  there is no assumption on the regularity
in variables $x$ and $y$. In this article,
we study  weighted $L^p$-norm inequalities  for spectral multipliers of $L$.
We show that sharp weighted H\"ormander-type
spectral multiplier theorems  follow  from Gaussian heat kernel bounds and appropriate  $L^2$ estimates of the kernels of the spectral multipliers.
These results are applicable to spectral multipliers for large classes of operators
 including
 Laplace operators acting on Lie groups of polynomial growth or irregular non-doubling domains
 of Euclidean spaces, elliptic
 operators on compact manifolds and Schr\"odinger operators with non-negative potentials.
\end{abstract}

\maketitle

 \tableofcontents

\bigskip

\section{Introduction}
\setcounter{equation}{0}

Suppose that $L$ is  a  non-negative self-adjoint operator  acting  on $L^2({X})$. Let $E(\lambda)$  be the spectral resolution
	  of  $L$. By the spectral theorem, for any bounded Borel function \linebreak $F: [0, \infty)\rightarrow {\Bbb C}$,  
one can define the operator
\begin{eqnarray}
F(L)=\int_0^{\infty} F(\lambda) dE(\lambda),
\label{e1.1}
\end{eqnarray}

 \noindent
 which is  bounded on $L^2(X)$.
A natural problem considered in the spectral multipliers theory is to give sufficient conditions on   $F$ and  $L$
  which imply the boundedness of $F(L)$  on various functional 
spaces defined on $X$. This topic has attracted a lot of attention and has been studied extensively by many authors: for example,
for sub-Laplacian on nilpotent groups in \cite{C}, \cite{D}, for sub-Laplacian on Lie groups of polynomial growth in \cite{A1},
for Schr\"odinger operator on Euclidean space $\mathbb R^n$ in \cite{He}, for sub-Laplacian on 
Heisenberg groups in \cite{MSt} and many others. For more information about the background of this topic,
the reader is referred to
 \cite{A1, A2, B,   C, DeM, D, DOS, FS, KW}  and the references therein. We also refer the reader to \cite{Str} and 
 the references therein for examples of potential applications of the spectral multiplier results. \\

We wish to point out \cite{DOS}, which is closely related to this paper. In \cite{DOS},
a sharp spectral  multiplier for a
non-negative self adjoint operator $L$ was obtained under  the assumption of the 
 kernel $p_t(x,y)$ of the analytic semigroup $e^{-tL}$ having 
a Gaussian upper bound. As there was no assumption on smoothness of the space variables
of $p_t(x,y)$, the singular integral $F(L)$ does not satisfy the standard kernel regularity condition
of a so-called Calder\'on-Zygmund operator, thus standard techniques of Calder\'on-Zugmund theory are not applicable. The lacking of smoothness of the kernel was indeed the main
obstacle in \cite{DOS} and it was overcome by shrewd exploitation of the analyticity of the
kernel $p_t(x,y)$ in variable $t$, together with a so-called Plancherel estimate, see Remark 2 after Corollary \ref{cor3.4}. \\

We will now recall  some of  main features of the spectral multipliers theory. 
An interesting example of a spectral multiplier  result comes from the paper \cite{A1} where Alexopoulos 
considers the operators acting on Lie groups of polynomial growth. He proved that 
if $L$ is a group invariant Laplacian and $n$ is the maximum of the local and global dimension of the group 
then  $F(L)$ is bounded on $L^p(\mathbb R^n)$ for all $1 < p < \infty$
if the function $F$ is differentiable $s$ times where $s = \left[\frac{n}{2}\right] + 1$ and satisfies
$$ 
 |\lambda^k F^{(k)}(\lambda) | \le C
 $$
for some constant $C$ and $k = 0, 1, \cdots , s$, see also Section~\ref{sec61} and Proposition~\ref{prop6.2} below. The philosophy is that we need function
$F$ to possess just more than  $n/2$ derivatives (with suitable bounds) for $F(L)$ to be bounded on all $L^p$ spaces,
$1 < p < \infty$. \\

When $s$ is an even number the above condition can be written in the following way
 \begin{equation}\label{Ale}
     \sup_{t>0}\| \eta \,\delta_t F \|_{W^\infty_s} < \infty,
     \end{equation}
   where $ \delta_t F(\lambda)=F(t\lambda)$,
   $\| F \|_{W^p_s}=\|(I-d^2/d x^2)^{s/2}F\|_{L^p}$ and $\eta$ is an auxiliary non-zero cut-off
    function such that $\eta \in C_c^{\infty}(\R_+)$. We note that condition (\ref{Ale}) is actually independent of the choice of $\eta$.
It is well known that condition (\ref{Ale}) can be generalized
with positive numbers $s>0$ and it is sufficient  
to take real value $s>n/2$, see \cite{A2, DOS}.  It is an interesting question when condition (\ref{Ale}) 
can be replaced
by the following weaker condition 
\begin{equation}\label{hor}
     \sup_{t>0}\| \eta \,\delta_t F \|_{W^2_s} < \infty
     \end{equation}
for some $s>n/2$.
Already for the standard Laplace operator on the Euclidean space $\mathbb R^n$, the classical Fourier multiplier result of 
H\"ormander \cite{Ho1} applied  to radial functions says that the weaker ${W^2_s}$ condition for any $s>n/2$ is enough to guarantee $L^p$ boundedness 
of $F(\Delta)$ for all $1<p < \infty$, see also \cite{C} for further discussion. 
 Actually, 
replacing the $W^\infty_s$ norm in condition  (\ref{Ale}) by the $W^2_s$ norm in condition (\ref{hor}) 
is essentially the same problem which one encounters in 
sharp Bochner-Riesz summability analysis, see \cite{CSo, DOS, Sog1, Sog2}. Discussion of possibility of replacing  condition 
(\ref{Ale}) by  (\ref{hor}) is one of the main themes of \cite{DOS}. \\

The aim of this  paper is to extend  the study  of sharp spectral multipliers in \cite{DOS} to the setting of weighted $L^p$ spaces. It turns out that for a function $F$ having more than 
$n/2$ suitable derivatives, the range of $p$ that we can obtain for $F(L)$ to be bounded depends also on the weight $w$. Most of the results of \cite{DOS} follow from 
Theorems \ref{th3.1}, \ref{th3.2} and \ref{th3.3}  which are the main results of this paper; see Remark 1 after Corollary \ref{cor3.4}. We use the techniques developed in \cite{DOS} to estimate the kernels of
spectral multipliers. The new contribution of this paper is a development of an original  technique to deal  with  
singular integral nature of the considered spectral multipliers to obtain generalization of 
unweighted results described in \cite{DOS} to weighted $L^p$ spaces.  \\

This paper was organized as follows. In Section 2, we recall basic properties of spaces of 
homogeneous type, the class of Muckenhoupt weights and a sufficient condition for boundedness of weighted singular integrals from \cite{AM}.  We state the  main results on 
weighted spectral multipliers, Theorems \ref{th3.1} and \ref{th3.3} in Section 3.
Section 4 is devoted to the proofs of these theorems. 
In Section 5, we use complex interpolation
to obtain boundedness for spectral multipliers on weighed $L^p$ spaces.
In Section 6, we give applications of our results to 
various operators in different settings, including Laplace operators on homogeneous groups
and on irregular domains of Euclidean spaces,
elliptic pseudo-differential operators on compact manifolds, Schr\"odinger operators with positive potentials 
 and holomorphic functional calculi of non-negative self-adjoint operators.\\


\section{Singular integrals and weights  }
\setcounter{equation}{0}

Let  $(X,d,\mu)$ be a space
endowed with a distance $d$ and a nonnegative Borel measure $\mu$ on $X$. Set $B(x,r)=\{y\in X: d(x,y)< r\}$ and  $V(x,r)=\mu(B(x,r))$.
 We shall  often just use $B$ instead of $B(x, r)$.
Recall that  $(X,d,\mu)$  satisfies the  doubling  volume property  provided that there exists a 
constant $C>0$ such that
\begin{eqnarray}
V(x,2r)\leq C V(x, r)\quad
\forall\,r>0,\,x\in X,
\label{e2.1}
\end{eqnarray}
%
%
more precisely if 
there exist $n, C_n>0$ such that
\begin{equation}
\frac{V(x,r)}{V(x,s)}\le C_n\left(\frac{r}{s}\right)^n, \quad 
\forall\,r\ge s>0,\,x\in X.
\label{e2.2}
\end{equation}
%
The parameter $n$ is a measure of the doubling  dimension of the space. It also follows from the doubling condition  that there  
exist $C$ and ${D}$, $0\leq {D}\leq n$ so that

\begin{equation}
V(y,r)\leq C\Big( 1+{\frac{d(x,y)}{ r}}\Big)^{D} V(x,r) \quad 
\forall\,r >0,\,x,y\in X.
\label{e2.3}
\end{equation}

\noindent
uniformly for all $x,y\in X$ and $r>0$. Indeed, property (\ref{e2.3}) with
${D}=n$ is a direct consequence of the triangle inequality for the metric
$d$ and  (\ref{e2.2}). In many cases like the Euclidean space
${\mathbb R}^n$ or Lie groups of polynomial growth, ${D}$ can be 
chosen to be $0$. 

\bigskip 

 \noindent {\bf Muckenhoupt weights.}
Next we review the definitions  of  Muckenhoupt classes of weights. We use the notation

$$
\oint_E h= {\frac{1}{ V(E)}} \int_E h(x)d\mu(x)
$$

\noindent
and we often forget the measure and variable of the integrand in writing integrals.

\medskip
 In what follow for any number or symbol $s$ with value in $[1, \infty]$ by $s'$ we denote it's conjugate, 
 that is $\displaystyle {\frac{1}{s}+\frac{1}{s'}=1}$.

A weight $w$  is a non-negative locally integrable function. We say that
$w\in A_p, 1<p<\infty$, if there exists a constant $C$ such that for every ball $B\subset X,$

$$
\Big( \oint_B  w \Big)  \Big( \oint_B  w^{1-p'} \Big)^{p-1} \leq C.
$$

 \noindent
For $p=1,$ we say that $w\in A_1$ if there is a constant $C$ such $ \Ma w\leq Cw$ a.e. 
where $\Ma$ denotes the uncentered maximal operator over balls 
  in $X$, that is
$$
\Ma w(x) =\sup_{B \ni x}\oint_B  w.
$$

 The reverse H\"older classes are defined in the following
way: $w\in RH_{q}, 1<q<\infty$, if there is a constant $C$ such that for every ball $B\subset X,$

$$
\Big( \oint_B w^q \Big)^{1/q}  \leq C\Big(  \oint_B  w \Big).
$$

\noindent
The endpoint $q=\infty$ is given by the condition: $w\in RH_{\infty}$ whenever, for any ball $B,$

$$
w(x)\leq C \oint_B  w, \ \ \ \ {\rm for \ a.e.}\  x\in B.
$$

\noindent
Note that we have excluded the case $q=1$ since the class $RH_1$ consists of all weights, and that is the way
$RH_1$ is understood in what follows.

We sum up some properties of the $A_p$ and $RH_q$ classes in the following lemmas.

\begin{lemma} \label{le2.1} Suppose that $(X,d,\mu)$ is a metric, measure space, which satisfies 
doubling condition {\rm (\ref{e2.1})}. Then the following properties hold for the weights classes 
$A_p$ and $RH_q$ defined on $(X,d, \mu)$: 
 
\smallskip
{\rm (i)}\ $A_1\subset A_p\subset A_q$ for   $1\leq p\leq q<\infty$.

\smallskip
{\rm (ii)}\ $RH_{\infty}\subset RH_q\subset RH_p$ for   $1\leq p\leq q\leq \infty$.

\smallskip
{\rm (iii)}\ If $w\in A_p, 1<p<\infty$, then there exists $1<q<p$ such that $w\in A_q$.

 \smallskip
{\rm (iv)}\ If $w\in RH_q, 1<q<\infty$, then there exists $q<p<\infty$ such that $w\in RH_p$.

\smallskip
{\rm (v)}\ $A_{\infty}=\cup_{1\leq p<\infty} A_p \subseteq \cup_{1<q\leq \infty} RH_q$.

 \smallskip
{\rm (vi)}\  If $1<p<\infty, w\in A_p$ if and only if $w^{1-p'}\in A_{p'}$.

 \smallskip
{\rm (vii)}\ If $1\leq q\leq \infty $ and $1\leq s<\infty$, then
$w\in A_q\cap RH_s$ if and only if $w^s\in A_{s(q-1)+1}$.
 \end{lemma}
 
\noindent
{\bf Proof}.  Properties (i)- (vi) are standard, see for instance, \cite{ST}, \cite{GR}  and \cite{Duo}.
For (vii), see \cite{JN}.
\hfill{}$\Box$

\medskip

Note that under   additional assumption on the measure $\mu$ that the function $\mu(B(x,r))$ increases   continuously with $r$
for each $x\in X$, it is shown that 
$A_{\infty}=\cup_{1\leq p<\infty} A_p = \cup_{1<q\leq \infty} RH_q$ (see  Theorem 18, Chapter 1, \cite{ST} ). 
However, we do not need this property
in the sequel.

\medskip

 \begin{lemma}\label{lemma2.2}
Let  $1<p<r^{'}$. Then  $w\in A_p\cap RH_{(\frac{{r^{'}}{ p}})^{'}}$
if and only if $w^{1-p'}=w^{-{\frac{1}{ p-1}}}\in A_{\frac{p'}{ r}}$. 
\end{lemma}

\medskip

\noindent
{\bf Proof}.  Lemma~\ref{lemma2.2} is a special case of \cite[Lemma 4.4]{AM} (with 
$p_0=1$ and $q'_0=r$ in the notation of \cite{AM}.)
\hfill{}$\Box$

\bigskip


\noindent 
{\bf Singular integrals on weighted spaces.}
The following result, see \cite[Theorem 3.7]{AM} is the main technical tool to
extend unweighted $L^p$ boundedness of spectral multipliers in 
\cite{DOS} to weighted $L^p$   results.

 \begin{theorem}\label{th2.4}
Let $1\leq p_0< \infty.$  Let $T$ be a sublinear operator acting on $L^{p_0}(X)$,
Let $\{A_r\}_{r>0}$ a family of operators acting on  $L^{p_0}(X)$. Assume that

\begin{eqnarray}\label{e2.4}
\Big( \oint_{B} \big| T(I-A_{r_B})f\big|^{p_0}d\mu\Big)^{1/p_0}
\leq C \Ma \big(\big| f\big|^{p_0}\big)^{\frac{1} {p_0}}(x)
\end{eqnarray}
and 
\begin{eqnarray}\label{e2.5}
 \big\| T A_{r_B}f \big\|_{L^\infty(B)}
\leq C \Ma \big(\big| Tf\big|^{p_0}\big)^{\frac{1} {p_0}}(x)
\end{eqnarray}

\noindent
for all $f \in L^{p_0}(X) $, and all ball $B$ with radius $r_B$ and all $B\ni x$.  
  Then for all 
 $p_0<p<\infty$  and $w\in A_{p/p_0}=A_{p/p_0}\cap RH_{ 1}$, 
there exists a constant $C$ such that
\begin{eqnarray}\label{e2.6}
   \|T f\|_{L^p(X, w)}\leq    C\|  f\|_{L^p(X, w)}.
\end{eqnarray}
\end{theorem}

\noindent
{\bf Proof}.  Theorem~\ref{th2.4} is a special case of \cite[Theorem 3.7]{AM} (with 
$q_0=\infty$ in the notation of \cite{AM}.)
\hfill{}$\Box$

\bigskip

Given $1\leq p_0<p<q_0$, we observe  that if $w$ is any given weight so that $w,
w^{1-p'}\in L^1_{\rm loc}(X),$ then a given linear operator $T$ is bounded on $L^p(X, w) 
$  if and only 
if its adjoint (with respect to $d\mu$) $T^{\ast}$ is bounded on $L^{p'}(w^{1-p'})$. Therefore,

\begin{eqnarray}\label{e2.7}
T: L^p(X, w)\rightarrow L^p(X, w)  \ \ \ {\rm for\ all\ } 
w\in A_{p\over p_0}\cap RH_{({q_0\over p})^{'}}
\end{eqnarray}

\noindent
if and only if 

\begin{eqnarray}\label{e2.8}
T^{\ast}: L^{p'}(X, w)\rightarrow L^{p'}(X, w)  \ \ \ {\rm for\ all\ } 
w\in A_{{p'}\over  q_0^{'}}\cap RH_{({p_0^{'}\over p'})^{'}}.
\end{eqnarray}

\medskip
  
The following result is a special case of interpolation with change  of measures. It was 
 proved in \cite{St} and \cite{SW} when $X={\Bbb R}^n$ is the Euclidean space.

 \begin{prop}\label{prop2.4}
Let  $1<r\leq q<\infty$ and let $w_0$ and $w_1$ be two positive weights. If $T$
is a bounded linear operator acting on  $L^r(X, w_0)$ and  $L^q(X, w_1)$.
Then $T$ is bounded on $L^p(X, w)$  for $r\leq p\leq q$
and $w=w_0^t w_1^{1-t}$, provided $t={q-p\over q-r}$ for $r\not=q$
and $0\leq t\leq 1$ for $r=q$.  
\end{prop}

\medskip

Note that $w^r\in A_p, r\geq 1$, if and only if $w\in A_p$
and $w$ satisfies $w\in RH_r$ and $w^{1-p'}=w^{-1/(p-1)}\in RH_r $ for  $p>1$; when $p=1$, we only need
$w\in RH_r$ 
(see pp. 351-352 of \cite{KW}).

\bigskip 
 
 \section{General spectral multiplier theorems on weighted spaces}
 \setcounter{equation}{0}

Let  $(X,d,\mu)$ be a space of homogeneous type. 
Recall that ${D}$ is the power that appeared in property (\ref{e2.3})
 and $n$ the  dimension entering  doubling volume condition   (\ref{e2.2}). 
 
Unless otherwise specified in the sequel we always assume that $L$ is a non-negative self-adjoint operator on $L^2(X)$
and that 
the semigroup $e^{-tL}$, generated by $-L$ on $L^2(X)$,  has the kernel  $p_t(x,y)$ 
which  satisfies 
the following  Gaussian upper bound

$$
\big|p_t(x,y)\big|\leq {C\over V(x,t^{1/m})} \exp\Big(-{  {d(x,y)^{m/(m-1)}\over  c\, t^{1/(m-1)}}}\Big)
\leqno{(GE)}
$$

 \noindent
for all $t>0$,  and $x,y\in X,$   where $C, c$ and $m$ are positive constants and $m\geq 2.$

\medskip

Such estimates are typical for elliptic or sub-elliptic differential operators of
order $m$ (see for instance, \cite{Da1}, \cite{DOS}, \cite{R} and \cite{VSC}).

  Theorems~\ref{th3.1}, \ref{th3.2} and \ref{th3.3} below  are the  main new results obtained in this paper. 
 
 \begin{theorem}\label{th3.1} \, Let $L$ be a non-negative self-adjoint operator such that the corresponding 
 heat kernels satisfy Gaussian bounds  $(GE)$.
 Let $s>{n\over 2}$ and let   $r_0=  \max\big(1, {2(n+{D})\over 2s+{D}}\big)$.
  Assume that for any $R>0$ and all Borel functions $F$  such that\, {\rm supp} $F\subseteq [0, R]$,

  \begin{eqnarray}\label{e3.2}
\int_X |K_{F(\sqrt[m]{L})}(x,y)|^2 d\mu(x) \leq {C\over V(y, R^{-1})} \|\delta_R F\|^2_{L^q} 
\end{eqnarray}

\noindent
for some $q\in [2, \infty].$   
 Then for any bounded Borel function $F$
 such that 
$\sup_{t>0}\|\eta\, \delta_tF\|_{W^q_s}<\infty,$
 the operator $F(L)$ is bounded on $L^p(X, w)$ for all $p$ and $w$ satisfying
  $  r_0  < p<\infty$ and   $w\in A_{p\over r_0}$. \noindent  
In addition,
 
\begin{eqnarray*} 
   \|F(L)  \|_{L^p(X, w)\rightarrow L^p(X, w)}\leq    C_s\Big(\sup_{t>0}\|\eta\, \delta_tF\|_{W^q_s}
   + |F(0)|\Big).
\end{eqnarray*}
\end{theorem}

\medskip

Note that Gaussian bounds $(GE)$ implies estimates (\ref{e3.2}) for $q=\infty$. This means that one can omit 
condition (\ref{e3.2}) if the case $q=\infty$ is consider. We describe the details in Theorem \ref{th3.2} below.

 \medskip
 
 \begin{theorem} \label{th3.2} \, Let $L$ be a non-negative self-adjoint operator such that the corresponding 
 heat kernels satisfy Gaussian bounds  $(GE)$. 
 Let $s>{n\over 2}$ and let   $r_0=  \max\big(1, {2(n+{D})\over 2s+{D}}\big)$.
 Then for any bounded Borel function $F$
 such that 
$\sup_{t>0}\|\eta\, \delta_tF\|_{W^{\infty}_s}<\infty,$
 the operator $F(L)$ is bounded on $L^p(X, w)$ for all $p$ and $w$ satisfying
  $  r_0  < p<\infty$ and   $w\in A_{p\over r_0}$. \noindent  
In addition,
 
\begin{eqnarray*} 
   \|F(L)  \|_{L^p(X, w)\rightarrow L^p(X, w)}\leq    C_s\Big(\sup_{t>0}\|\eta\, \delta_tF\|_{W^{\infty}_s}
   + |F(0)|\Big).
\end{eqnarray*}
 \end{theorem}
 
 \bigskip

  \noindent
{\bf Proof.}  Note that it was proved in Lemma~2.2 of \cite{DOS}, that
for any Borel function $F$  such that supp $F\subset [0, R],$

\begin{eqnarray}\label{eff}
\big\|K_{F(\sqrt[m]{L})}(\cdot, y)\big\|^2_{L^2(X)} &=& 
\big\|K_{{\overline{F}}(\sqrt[m]{L})}(y, \cdot)\big\|^2_{L^2(X)}\nonumber\\
&\leq& {C \over V(y, R^{-1}) }
\big\|F\big\|^2_{L^{\infty}}
\end{eqnarray}
where ${\overline{F}}$ denotes the complex conjugate of $F$. 

 This shows that estimate (\ref{e3.2}) always holds for $q=\infty,$ and Theorem~\ref{th3.2}  follows from
 Theorem~\ref{th3.1}.
 \hfill{} $\Box$
 
 \medskip

  From a point of view of some applications of spectral multipliers the sharp results and  the  required number of derivatives are not essential for the final outcome, see for example \cite{Str}. For this kind  of 
 applications Theorem~\ref{th3.2} is the  best solution because to use it one does not have to consider or prove 
 condition (\ref{e3.2}). 
 Nevertheless, Theorem~\ref{th3.1} and condition (\ref{e3.2}) is of  significant interest independent of their  applications.  
 In the case of standard Laplace operator condition (\ref{e3.2}) is equivalent with $(1,2)$ restriction theorem
 and both Theorem~\ref{th3.1} and condition (\ref{e3.2}) are a new part of Bochner-Riesz analysis. 
 Estimates (\ref{e3.2}) are also closely related to Strichartz and other dissipative type  estimates. 
 For further discussion of condition (\ref{e3.2}), see also \cite{DOS}.

It is not difficult to see  that condition  (\ref{e3.2}) with some $q< \infty$ implies that the set of point spectrum of the considered operator is empty  because 
the $L^q$ norm of characteristic function of any singleton subset of $\R$ is zero. Hence  if $q<\infty$ then $F(\sqrt[m]{L})$ does not depend on the value of $F(0)$ because then the point spectrum is empty and the spectral projection 
on zero eigenvalue $E(\{0\})=0$. 
Therefore if $q < \infty$ then one can skip $|F(0)|$ in the concluding estimates of Theorem~\ref{th3.1}. 
See \cite[(3.3)]{DOS} for more detailed explanation.   

The fact that the  point spectrum of the considered operator is empty implies also that for elliptic operators on compact manifolds condition  (\ref{e3.2}) cannot hold for any $q<\infty$.
To be able to  study these operators as well, similarly as in  \cite{CS, DOS} we introduce some variation of condition  (\ref{e3.2}).
For a Borel function $F$ such that supp $F\subseteq [-1, 2]$ we define
the norm $\|F\|_{N,q}$ by the formula

$$
\|F\|_{N,q}=\Big({1\over 3N}\sum_{\ell=1-N}^{2N} \sup_{\lambda\in
 [{\ell-1\over N}, {\ell\over N})} |F(\lambda)|^q\Big)^{1/q},
 $$

 \noindent
 where $q\in [1, \infty)$ and $N\in{\mathbb Z}_+$. For $q=\infty$, we put 
 $\|F\|_{N, \infty}=\|F\|_{L^{\infty}}$.
 It is obvious that $\|F\|_{N,q}$ increases monotonically in $q$. 
 
 The next theorem
 is a variation of Theorem~\ref{th3.1}. This variation can be used in case of operators with nonempty
 point spectrum, see also \cite[Theorem 3.6]{CS} and \cite[Theorem 3.2]{DOS}.

 \medskip

 \begin{theorem}\label{th3.3} 
 Assume that $\mu(X)<\infty$.  \, Let $L$ be a non-negative self-adjoint operator such that the corresponding 
 heat kernels satisfy Gaussian bounds $(GE)$. 
 Let $s>{n\over 2}$ and let   $r_0=  \max\big(1, {2(n+{D})\over 2s+{D}}\big)$.
Suppose  
  that for any $N\in{\mathbb Z}_+$ and for all  
  Borel functions $F$  such that\, {\rm supp} $F\subseteq [-1, N+1]$, 
  
   \begin{eqnarray}\label{e3.4}
\int_X |K_{F(\sqrt[m]{L})}(x,y)|^2 d\mu(x) \leq {C\over V(y, N^{-1})} \|\delta_N F\|^2_{N,\, q} 
\end{eqnarray}

\noindent
for some $q\geq 2.$ 
  Then for any bounded Borel function $F$ such that $\sup_{t>1}\|\eta\, \delta_tF\|_{W^q_s}<\infty,$
  the operator $F(L)$ is bounded on $L^p(X, w)$ for all $p$ and $w$ satisfying
  $  r_0  < p<\infty$ and   $w\in A_{p\over r_0}$. \noindent  
In addition,
 
\begin{eqnarray*} 
   \|F(L)  \|_{L^p(X, w)\rightarrow L^p(X, w)}\leq    C_s\Big(\sup_{t>1}\|\eta\, \delta_tF\|_{W^q_s}
   +\|F\|_{L^{\infty}}\Big).
\end{eqnarray*}
 \end{theorem}

 \medskip
 
 We will discuss  the proofs of Theorems \ref{th3.1} and \ref{th3.3} in Section 4. These results have the following 
 corollary.
 
 \medskip
 
 \begin{cor} \label{cor3.4}
 Let $s>{n\over 2}$ and let   $r_0=  \max\big(1, {2(n+{D})\over 2s+{D}}\big)$ 
 and ${1\over r_0} +{1\over r'_0}=1$.
 Suppose in addition  that $1<p< r_0^{'}$ and   
 $w\in A_{p}\cap RH_{ ({r_0^{'}\over p} )^{'}}$.
 
 \noindent
 {\rm (a)} Assume also that the operator $L$ satisfies the assumptions of Theorem {\rm \ref{th3.1}}
  for some $2 \le q \le \infty$, then
  \begin{eqnarray*} 
    \|F(L)  \|_{L^p(X, w)\rightarrow L^p(X, w)}\leq    C_s\Big(\sup_{t>0}\|\eta\, \delta_tF\|_{W^q_s}
    + |F(0)|\Big).
 \end{eqnarray*}
{\rm  (b)} Alternatively assume in addition that  the operator $L$ satisfies the assumptions of Theorem  {\rm \ref{th3.3}} 
 for some $2 \le q \le \infty$, then 
 \begin{eqnarray*} 
    \|F(L)  \|_{L^p(X, w)\rightarrow L^p(X, w)}\leq    C_s\Big(\sup_{t>1}\|\eta\, \delta_tF\|_{W^q_s}
    +  \|F\|_{L^{\infty}}\Big).
    \end{eqnarray*}
 \end{cor}

\smallskip

 \bigskip

  \noindent
{\bf Proof.} Suppose  $1<p< r_0^{'}$ and   
 $w\in A_{p}\cap RH_{ ({r_0^{'}\over p} )^{'}}$.   We have that $w^{-{1\over p-1}}\in A_{p^{'} \over r_0}.$
Then for $f\in L^{\infty}_c(X)$ ({\it i.e.} bounded with compact
  support),
  
    \begin{eqnarray*}
  \|F(L)f\|_{L^p(X, w)}
  &=& \Big|\int_X  F(L)f(x){\overline{  g(x)}} d\mu(x)\Big|,
  \end{eqnarray*}
  
  \noindent
  where the supremum is taken over all functions $g\in L^{\infty}_c(X)$  
   such that  $\|g\|_{L^{p'}(X,w^{-{1\over p-1}})}=1.$
  
  Let ${\bar{ F}(L) } $ be the operator with multiplier ${\bar{ F}} $, the complex conjugate
  of $F$. Then  ${\bar { F}} $ satisfies the same estimates as $F$, and we have 
  \begin{eqnarray*}
  \|F(L)f\|_{L^p(X, w)}
  &= &\sup 
  \Big|\int_X  f(x) {\bar{ F}(L) }\,   {{ g(x)}} d\mu(x)\Big|\nonumber\\[2pt]
    &\leq& \sup \|f\|_{L^{p}(X,w)}
  \big\|  {\bar{ F}(L) }\,g\|_{L^{p'}(X, w^{-{1\over p-1}} )}
    \nonumber\\[2pt]
    &\leq & C\|f\|_{L^{p}(X,w)}
  \end{eqnarray*}
  
  \noindent
    since   $p'> r_0$, and we can apply Theorems \ref{th3.1} or  \ref{th3.3} to the weight
  $w^{-{1\over p-1}}\in A_{p^{'} \over r_0}.$
 \hfill{} $\Box$
\bigskip

\noindent
  {\bf Remarks.}  
  
\smallskip

1)  \, Note that Theorems~\ref{th3.1} and  \ref{th3.3} imply the main results obtained in  \cite{DOS}.
Indeed the  trivial weight $w=1$ is in all $A_p$ classes, so under the assumptions  of Theorems~\ref{th3.1} and  \ref{th3.3}
the operator $F(L)$ is bounded on all $L^p$ spaces  $1<p<\infty$. Note that for $p<2$, 
$L^p$ boundedness of $F(L)$ follows by considering the adjoint operator $F(L)^*=\bar{F}(L)$.  
 Similarly to the results in \cite{DOS} the important point of this paper is that if one can obtain (\ref{e3.2}) or (\ref{e3.4})
 for some $q < \infty$ then one can prove stronger multiplier results than in case $q=\infty$. The estimates (\ref{e3.2})
 for $q=\infty$ are not necessary because estimates   (\ref{e3.2}) with $q=\infty$ follow from  Gaussian bounds assumption $(GE)$, 
 see Theorem \ref{th3.2} . If one has (\ref{e3.2}) or (\ref{e3.4}) for $q=2$,
then this implies the sharp weighted  H\"ormander-type multiplier result. Actually,
we believe that to obtain any sharp weighted  H\"ormander-type multiplier theorem one has
to investigate conditions of the same type as (\ref{e3.2}) or (\ref{e3.4}), i.e.
conditions which allow us to estimate the norm $\big\|K_{F{\sqrt[m]{L}}}(\cdot, y)\big\|^2_{L^2(X, \mu)}  
$ in terms of some kind of $L^p$ norm of the function $F.$

\medskip

2) \,We call hypothesis (\ref{e3.2}) or (\ref{e3.4}) the Plancherel estimates or the Plancherel
conditions. For the standard Laplace operator on Euclidean spaces $\mathbb R^n$, this is equivalent to 
$(1,2)$ Stein-Tomas restriction theorem (which is also the Plancherel estimate of the Fourier transform).   
Assumption that $q \ge 2$ is not necessary in the proofs of  Theorems~\ref{th3.1}  and ~\ref{th3.3}. 
However we do not expect that there are any examples where  estimates (\ref{e3.2}) or (\ref{e3.4})
hold with $q<2$ because this would imply the Riesz summability for the index $\alpha<(n-1)/2$
which is  false for the standard Laplace operator. 

\medskip



3)  If we take $s>n/2$ in Theorems~\ref{th3.1} and ~\ref{th3.3}, 
 then for every $w\in A_1\cap RH_2$, the operator $F(L)$ maps $L^1(X, w)$ into $L^{1, \infty}(X, w)$, that is,
 there is a constant $C>0$, independent of $f$ and $\lambda$, such that 
 \begin{eqnarray*}
 w\big\{x\in X: \,  \big|F(L)f(x)\big|>\lambda \big\}\leq {C\over \lambda} \big\|f\big\|_{L^1(X, w)},
 \ \ \ \ \ \lambda>0.
 \end{eqnarray*} 
 
 \noindent
 The proof follows from the line  of Theorem 5.8 in \cite{M}, together
with the proofs of Theorems 3.1 and 3.2 in \cite{DOS}, respectively. The details are left to the reader.


\bigskip

\section{Proofs  of Theorems~\ref{th3.1} and \ref{th3.3}}
\setcounter{equation}{0}

Recall that  $B=B(x_B, r_B)$ is the ball of radius $ r_B$ at centred at $ x_B$. Given 
$\lambda>0$, we will write $\lambda B$  the ball with the same centre as $B$ and with radius $r_{\lambda B}=\lambda r_B$. 
We set

\begin{equation}
U_0(B)=B,  \ \  {\rm and}\ \ \ U_j(B)=2^j B\backslash 2^{j-1}B 
\,\,\mbox{ for }\,\,j=1,2, \dots.  
\label{e4.1}
\end{equation}

\smallskip

As a preamble to the proof of Theorem~\ref{th3.1},  we record
a useful auxiliary result. For a proof, see pp. 453-454, Lemma 4.3 of \cite{DOS}.

\begin{lemma}\label{le4.1} {\rm (a)}\, Suppose that $L$ satisfies {\em (\ref{e3.2})} for some
$q\in [2, \infty]$ and that $R>0, s>0$. Then for any $\epsilon>0$, there exists a constant
$C=C(s, \epsilon)$ such that 
 
\begin{eqnarray}\label{e4.2}
\int_X \big|K_{F(\sqrt[m]{L})}(x,y)\big|^2 \big(1+Rd(x,y)\big)^{s} d\mu(x)\leq 
{C \over V(y, R^{-1})}
 \|\delta_{R} F\|^2_{W^q_{{s\over 2} +\epsilon}} 
\end{eqnarray}

\noindent
for all Borel functions $F$ such that supp $F\subseteq [R/4, R].$

  {\rm (b)}\, Suppose that $L$ satisfies {\rm (\ref{e3.4})} for some
 $q\in [2, \infty]$ and that $N>8$ is a natural number.
  Then for any $s>0,$ $\epsilon>0$
 and function $\xi\in C_c^{\infty}([-1, 1])$  there exists a constant
 $C=C(s, \epsilon, \xi)$ such that 
 
 \begin{eqnarray} \label{e4.3}
 \int_X \big|K_{F\ast \xi\, (\sqrt[m]{L})}(x,y)\big|^2 \big(1+Nd(x,y)\big)^{s} d\mu(x)\leq 
 {C \over V(y, N^{-1})}
 \|\delta_{N} F\|^2_{W^q_{{s\over 2} +\epsilon}} 
 \end{eqnarray}

 \noindent
 for all Borel functions $F$ such that supp $F\subseteq [N/4, N].$
\end{lemma}

\medskip

\noindent
 {\bf Proof   of Theorem~\ref{th3.1}.  }
 We fix  $s$ such that 
 $s>{\frac{n}{2}}$, and thus  $  {2(n+{D})\over 2s+{D}}<2.$  
 In this case, we take one parameter  $p_0$   in the sequel
such that  $p_0$ belongs to the interval $ \left(\max \left\{{2(n+{D})\over 2s+{D}}, 1\right\},  2\right)$. 
  Let $M\in{\mathbb N}$ such that $M>s/m,$  where $m$ is the constant in $(GE)$.  
 We will show that  
 for all balls $B\ni x$,      

\begin{eqnarray}
 \Big(  \oint_{B}
\big|F(L)(I-e^{-r_B^mL})^M f \big|^{p_0}d\mu\Big)^{1/p_0} 
\leq C  \Ma \big(  \big| f\big|^{p_0}\big)^{1\over p_0}(x)
\label{e4.4}
\end{eqnarray}

\noindent
  for all $f\in L^{\infty}_c(X)$.

Let us prove (\ref{e4.4}). 
Observe that  $\sup_{t>0}\|\eta\, \delta_tF\|_{W^p_s} \sim \sup_{t>0}\|\eta\, \delta_t G\|_{W^p_s}$
where $G(\lambda)=F(\sqrt[m]{\lambda})$.
 For
this reason, we can replace $F(L)$ by  $F(\sqrt[m]{L})$ in the proof.
Notice that  $F(\lambda)=F(\lambda)-F(0)+F(0)$ and hence
   
   $$
   F(\sqrt[m]{L})=(F(\cdot)-F(0))(\sqrt[m]{L}) +F(0)I.
   $$
   
   \noindent
 Replacing $F$ by $F-F(0)$, we may assume in the sequel that $F(0)=0.$ 
   Let   $\varphi\in C_c^{\infty}(0, \infty)$ be a non-negative function satisfying 
    supp $\varphi \subseteq [{1\over 4}, 1]$ and
$
\sum_{\ell=-\infty}^{\infty}\varphi(2^{-\ell}\lambda) =1$ for any $\lambda >0,
$ 
and let $\varphi_\ell$ denote the function $\varphi(2^{-\ell}\cdot).$
Then

\begin{eqnarray}
F(\lambda)=
\sum_{\ell=-\infty}^{\infty}\varphi(2^{-\ell}\lambda)F(\lambda)
=  \sum_{\ell=-\infty}^{\infty}F^{\ell}(\lambda), \ \ \  \ \forall\,  \lambda \ge 0.
\label{e4.5}
\end{eqnarray}

\noindent
This decomposition implies that 
  the sequence $\sum_{\ell=-N}^{N}F^{\ell}(\sqrt[m]{L})$
converges strongly in $L^2(X)$ to $F(\sqrt[m]{L})$ (see for instance, Reed and Simon  \cite{RS}, Theorem VIII.5). 
  For every   $\ell\in{\Bbb Z}$ and $r>0$, we set for $\lambda>0,$

\begin{eqnarray}
F_{r, M}(\lambda)&=& F(\lambda)(1-e^{-(r\lambda)^m})^{M}\label{e4.6},\\[5pt]
F^{\ell}_{r, M}(\lambda)&=& F^{\ell}(\lambda)(1-e^{-(r\lambda)^m})^{M}.\label{e4.7}
\end{eqnarray}

\noindent
Given a  ball $B\subset X$, we  use the 
 decomposition
$
 f=\sum\limits_{j=0}^{\infty} f_j  $ in which $ f_j=f\chi_{U_j (B)},
 $
and $U_j(B)$ were defined in (\ref{e4.1}).  
We may  write
\begin{eqnarray}
\label{e4.8} 
 F(\sqrt[m]{L})(1-e^{- r_B^mL})^{M}f &=& F_{r_B, M}(\sqrt[m]{L})f \\
 &= &
 \sum_{j=1}^2 F_{r_B, M}(\sqrt[m]{L})f_j +    
 \lim_{N\rightarrow \infty}\sum_{\ell=-N}^{N}  \sum_{j= 3}^{\infty}  
F^{\ell}_{r_B, M}(\sqrt[m]{L})f_j,   \nonumber
\end{eqnarray}

\noindent
where the sequence  
converges strongly in $L^2(X)$. 

From Gaussian condition $(GE)$, we have  that for any $t>0$, 
$\|e^{-tL}f\|_{L^p(X)}\leq 
C\| f\|_{L^p(X)}$.
 This, 
 in combination with  $L^{p}$-boundedness of the operator $F(\sqrt[m]{L})$
 (see Theorem 3.1, \cite{DOS}), gives  that  
 for all balls $B\ni x$,      
\begin{eqnarray}\label{e4.9}
 \Big(  \oint_{B}
\big|F_{r_B, M}(\sqrt[m]{L})f_j \big|^{p_0}d\mu\Big)^{1/p_0} 
&\leq& V(B)^{-{1/p_0}} \big\|F_{r_B, M}(\sqrt[m]{L})f_j\big\|_{L^{p_0}(X)}\nonumber\\
&\leq&  CV(B)^{-1/p_0} \big\|  f_j\big\|_{L^{p_0}(X)}\nonumber\\
&\leq& C  \Ma \big(  \big| f\big|^{p_0}\big)^{1\over p_0}(x) 
\end{eqnarray}

\noindent
  for  $j=1,2.$

  Fix $j\geq 3.$ Let $p_1\geq 2$ and ${1\over p_0}-{1\over p_1}={1\over 2}.$
 By H\"older's inequality,  
 we have  that for all balls $B\ni x, $ 
 \begin{eqnarray}
\label{e4.10} \hspace{1cm}
 &&\hspace{-2cm} \Big(  \oint_{B}
\big|F^{\ell}_{r_B, M}(\sqrt[m]{L})  f_j \big|^{p_0}d\mu\Big)^{1/p_0} \\ [2pt] 
 &\leq& V(B)^{ -{1\over p_1}}  \big\|F^{\ell}_{r_B, M}(\sqrt[m]{L})  f_j\big\|_{L^{p_1}(B)}
\nonumber\\ [2pt]
  &\leq& V(B)^{ -{1\over p_1}}  \big\|F^{\ell}_{r_B, M}(\sqrt[m]{L})  \big\|_{L^{p_0}(U_j(B))\to  L^{p_1}(B)}
 \big\|  f_j\big\|_{L^{p_0}(X)}
 \nonumber\\ [2pt]
  &\leq& C2^{jn\over p_0 }   V(B)^{{1\over 2} } \big\|F^{\ell}_{r_B, M}(\sqrt[m]{L})  \big\|_{L^{p_0}(U_j(B))\to  L^{p_1}(B)}
  \Ma \big(  \big| f\big|^{p_0}\big)^{1\over p_0}(x). 
 \nonumber 
\end{eqnarray}

\noindent
Let  ${1\over p_0}={\theta\over 1} +{1-\theta\over 2}$ and ${1\over p_1}={\theta\over 2},$   that is $\theta=2({1\over p_0}-{1\over 2})$.
By interpolation, 

 \begin{eqnarray}
\label{e4.11} \hspace{1cm}
 &&\hspace{-2cm} \big\|{F}^{\ell}_{r_B, M}(\sqrt[m]{L})  \big\|_{L^{p_0}(U_j(B))\to  L^{p_1}(B)} \nonumber\\ [3pt]
 &\leq&      \big\|F^{\ell}_{r_B, M}(\sqrt[m]{L})  \big\|^{1-\theta}_{L^{2}(U_j(B))\to  L^{\infty}(B)}
    \big\|\overline{F}^{\ell}_{r_B, M}(\sqrt[m]{L})  \big\|^{\theta}_{L^{2}(B)\to  L^{\infty}(U_j(B))}.
\end{eqnarray}

\noindent
Next we estimate $\big\|F^{\ell}_{r_B, M}(\sqrt[m]{L})  \big\|_{L^{2}(U_j(B))\to  L^{\infty}(B)}.$
For every   $\ell\in{\mathbb Z}$,  let $K_{F^{\ell}_{r_B, M}(\sqrt[m]{L})}(y,z)$ be  the Schwartz  kernel of operator 
 $F^{\ell}_{r_B, M}(\sqrt[m]{L})$. Then we have

\begin{eqnarray}
  \hspace{-1.3cm}&&\hspace{-1.5cm} \big\|F^{\ell}_{r_B, M}(\sqrt[m]{L})  \big\|^2_{L^{2}(U_j(B))\to  L^{\infty}(B)} \nonumber\\ [3pt]\label{e4.12}
  &=& \sup_{y\in B}  
\int_{U_j(B)}  \big| K_{F^{\ell}_{r_B, M}(\sqrt[m]{L})}(y,z) \big|^2
   d\mu(z) \\ [3pt]
   &\leq& C  
 2^{-2s j} \big(2^{\ell} r_B\big)^{-2s } \sup_{y\in B} \int_{X} \big|K_{F^{\ell}_{r_B, M}(\sqrt[m]{L})}(y,z) 
  \big|^2 \big(1+2^{\ell}d(y, z)\big)^{2s }d\mu(z).
 \nonumber 
\end{eqnarray}

\noindent
We then apply Lemma~\ref{le4.1}   with 
 $F=F^{\ell}_{r_B, M} $ and $R=2^{\ell}$ to obtain 
\begin{eqnarray} 
  \hspace{1.5cm}
 \int_{X} 
  \big|K_{F^{\ell}_{r_B, M}(\sqrt[m]{L})} (y,z)
 \big|^2  \big(1+2^{\ell} d(y,z)\big)^{2s}d\mu(z) 
  \leq  {C_{s} \over
 V(y, 2^{-\ell})}\,
 \|\delta_{2^{\ell}} \big(F^{\ell}_{r_B, M}\big)\|^{2}_{W^q_s}.   
 \label{e4.13}
\end{eqnarray}
 
 \noindent
 Now for any Sobolev space $W^q_s({\Bbb R})$, if $k$ is an integer greater than $s$, then

\begin{eqnarray}\label{e4.14}
 \|\delta_{2^{\ell}}\big(F^{\ell}_{r_B, M}\big)\|_{W^q_s}
&=&\big\|\varphi(t)F(2^{\ell}t) (1-e^{-(2^{\ell}r_Bt)^m})^{M}\big\|_{W^q_s}  \nonumber\\ [4pt] 
  &\leq&   C \big\|(1-e^{-(2^{\ell} r_B t)^m})^M\big\|_{C^k([{1\over 4}, 1])} \,
  \big\|\delta_{2^{\ell}}[\varphi_{\ell} F]\big\|_{W^q_s}\nonumber\\ [4pt]
    &\leq&   
 C {\rm min}\, \big\{1, (2^{\ell}r_B)^{mM}\big\} \big\|\delta_{2^{\ell}}[\varphi_{\ell} F]\big\|_{W^q_s}.
\end{eqnarray}

\noindent
Note that for all $y\in B$, $B \subset B(y,2r_B)$ so by  (\ref{e2.2}) 
\begin{eqnarray}\label{e4.15}
 {1\over
 V(y, 2^{-\ell})} 
\leq   C\sup_{y\in B} {V(y, 2r_B)\over 
 V(y, 2^{-\ell}) V(B)  }  
   \leq   {C\over V(B)}\max\big\{1, \big(2^{\ell} r_B\big)^n\big\}.
\end{eqnarray}
Hence by  (\ref{e4.14}) and   (\ref{e4.15}), 

\begin{eqnarray}\label{e4.16}
 &&\big\|F^{\ell}_{r_B, M}(\sqrt[m]{L})  \big\|_{L^{2}(U_j(B))\to  L^{\infty}(B)}
 \\[3pt]
&\leq &C  \left( 
 2^{-2sj} \big(2^{\ell} r_B\big)^{-2s}   {\rm min}\, \big\{1, (2^{\ell}r_B)^{2mM}\big\}
 \max\big\{1, \big(2^{\ell} r_B\big)^n\big\} {1\over V(B)} \right)^{{1/ 2}}
 \big\|\delta_{2^{\ell}}[\varphi_{\ell} F]\big\|_{W^q_s}. \nonumber 
\end{eqnarray}

\noindent
We now turn  to  estimate the term $ \big\|\overline{F}^{\ell}_{r_B, M}(\sqrt[m]{L})  \big\|_{L^{2}(B)\to  L^{\infty}(U_j(B))}.$ 
The calculations symmetric to (\ref{e4.12}), (\ref{e4.13})  and  (\ref{e4.14}) with $\sup_{y\in B}$ replaced 
by $\sup_{z\in U_j(B)}$
 yields,  

\begin{eqnarray*}
 &&\hspace{-1cm}  \big\|\overline{F}^{\ell}_{r_B, M}(\sqrt[m]{L})  \big\|_{L^{2}(B)\to  L^{\infty}(U_j(B))}\\[2pt] 
  &\leq&   C\left(  2^{-2sj} \big(2^{\ell} r_B\big)^{-2s} {\rm min}\, \big\{1, (2^{\ell}r_B)^{2mM}\big\}
    \sup_{z\in U_j(B)} {1\over
 V(z, 2^{-\ell})}
 \right)^{1/2}\big\|\delta_{2^{\ell}}[\varphi_{\ell} F]\big\|_{W^q_s}.
\end{eqnarray*}
Next by (\ref{e2.2}) and (\ref{e2.3})
\noindent
\begin{eqnarray*}
 \sup_{z\in U_j(B)} {1\over
 V(z, 2^{-\ell})} 
&\leq &  C\sup_{z\in U_j(B)} \Bigg({V(z, r_B)\over 
 V(z, 2^{-\ell})}   \times   \Big( 1+{d(z, x_B)\over r_B}\Big)^{D}\Bigg) {1\over
 V(x_B, r_B)}\nonumber\\ [3pt]
   &\leq &  C{2^{j{D}}\over V(B)}\max\big\{1, \big(2^{\ell} r_B\big)^n\big\}.
\end{eqnarray*}
Hence 
\begin{eqnarray}\label{e4.17}
&&  \hspace{1cm}    \big\|\overline{F}^{\ell}_{r_B, M}(\sqrt[m]{L})  \big\|_{L^{2}(B)\to  L^{\infty}(U_j(B))}\\[3pt]
&\leq& C \left( 
 2^{-2sj} \big(2^{\ell} r_B\big)^{-2s} 2^{j{D}}  {\rm min}\, \big\{1, (2^{\ell}r_B)^{2mM}\big\}
 \max\big\{1, \big(2^{\ell} r_B\big)^n\big\} {1\over V(B)}  \right)^{ { 1/2}}
 \big\|\delta_{2^{\ell}}[\varphi_{\ell} F]\big\|_{W^q_s}. \nonumber
\end{eqnarray}

\noindent
It then follows from estimates  (\ref{e4.16}) and  (\ref{e4.17}), in combination with (\ref{e4.11}) and (\ref{e4.10}) that

\begin{eqnarray}\label{e4.18}
  &&\hspace{-2cm}  \Big(  \oint_{B}
\big|F^{\ell}_{r, M}(\sqrt[m]{L}) f_j \big|^{p_0}d\mu\Big)^{1/p_0} 
\nonumber\\[2pt]
 &\leq&  C  
2^{-js +{jn\over p_0} + {j{D} \theta\over 2 }} \left(
\big(2^{\ell}r_B \big)^{-s}  {\rm min}\, \big\{1, (2^{\ell}r_B)^{mM}\big\}
 \max\big\{1, \big(2^{\ell} r_B\big)^{  n\over 2}\big\}  \right) \nonumber\\[2pt]
   && \hspace{5cm} \times \Ma \big(\big|f\big|^{p_0}\big)^{1\over p_0}(x)
\sup_{\ell\in{\Bbb Z}}\big\|\delta_{2^{\ell}}[\varphi_{\ell} F]\big\|_{W^q_s}.
\end{eqnarray}

 \noindent
 Therefore,
 
 \begin{eqnarray} \label{e4.19}
 &&\hspace{-1.2cm}  
\sum_{j= 3}^{\infty}  
\sum_{\ell=-\infty}^{\infty}
\Big(  \oint_{B}
\big|F^{\ell}_{r, M}(\sqrt[m]{L}) f_j \big|^{p_0}d\mu\Big)^{1/p_0} 
\nonumber\\[2pt]
\hspace{1.2cm} &\leq&  C 
\sum_{j= 3}^{\infty}   
2^{-js +{jn\over p_0} + {j{D} \theta\over 2 }}\Big(\sum_{\ell=-\infty}^{\infty} 
\big(2^{\ell}r_B \big)^{-s}  {\rm min}\, \big\{1, (2^{\ell}r_B)^{mM}\big\}
 \max\big\{1, \big(2^{\ell} r_B\big)^{  n\over 2}\big\}    \Big)  \nonumber\\[2pt]
   && \hspace{1cm}\times \Ma \big(\big|f\big|^{p_0}\big)^{1\over p_0}(x)
\sup_{\ell\in{\Bbb Z}}\big\|\delta_{2^{\ell}}[\varphi_{\ell} F]\big\|_{W^q_s}\nonumber\\[2pt]
&\leq&  C 
\sum_{j= 3}^{\infty}   
2^{({n+{D}\over p_0}- (s+{{D}\over 2}))j} \Big(\sum_{\ell:\, 2^{\ell}r_B> 1  }\big(2^{\ell} r_B\big)^{-s+{n\over 2}}
 + \sum_{\ell:\, 2^{\ell}r_B\leq 1  }\big(2^{\ell}r_B\big)^{mM-s}   \Big)  \nonumber\\[2pt]
   && \hspace{1cm}\times \Ma \big(\big|f\big|^{p_0}\big)^{1\over p_0}(x)
\sup_{\ell\in{\Bbb Z}}\big\|\delta_{2^{\ell}}[\varphi_{\ell} F]\big\|_{W^q_s}\\[2pt]
&\leq&  C 
\sum_{j= 3}^{\infty}   
2^{({n+{D}\over p_0}- (s+{{D}\over 2}))j} 
   \Ma \big(\big|f\big|^{p_0}\big)^{1\over p_0}(x)
\sup_{\ell\in{\Bbb Z}}\big\|\delta_{2^{\ell}}[\varphi_{\ell} F]\big\|_{W^q_s}\nonumber\\[2pt]
 &\leq& C   
 \Ma \big(\big|f\big|^{p_0}\big)^{1\over p_0}(x)\sup_{\ell\in{\Bbb Z}}\big\|\delta_{2^{\ell}}[\varphi_{\ell} F]\big\|_{W^q_s}.\nonumber
\end{eqnarray}

\noindent
   Here, the second inequality is obtained by using condition  $\theta=2({1\over p_0}-{1\over 2})$,
   and the third inequality  follows from the convergence of power series with common ratio $1/2.$  
In the last inequality we have used the fact that   
  $p_0>{2(n+{D})\over 2s+{D}}$. 
 
 \medskip
 
Combining estimates  (\ref{e4.9}) and (\ref{e4.19}), we have therefore proved  (\ref{e4.4}),
and then estimate (\ref{e2.4}) holds for $T=F(L)$ and $A_{r_B}=I-(I-e^{-r_B^mL})^M$.
Note also that estimate (\ref{e2.5}) always holds for $A_{r_B}=I-(I-e^{-r_B^mL})^M$.
Indeed note that $T=F(L)$ and $A_{r_B}=I-(I-e^{-r_B^mL})^M$ commutes so it is enough to show that
$$
 \big\| A_{r_B}f \big\|_{L^\infty(B)}
\leq C \Ma \big(\big| f\big|^{p_0}\big)^{1\over p_0}(x).
$$
It is not difficult to see that it is enough to prove the above inequality for $p_0=1$. However 
 $A_{r_B}=I-(I-e^{-r_B^mL})^M$ is a finite linear combination of the terms $e^{-jr_B^mL}$, $j=1,\ldots M$
 which all satisfy Gaussian bounds and the above inequality and in turn (\ref{e2.5}) follow from that 
 observation. 

It then follows from Theorem~\ref{th2.4}   that
for all $p>p_0>r_0={2(n+{D})\over 2s+{D}}$, the operator $F(L)$ is bounded on $L^p(X, w)$ provided that
$ 
w\in A_{p\over p_0}.
$ 
On the other hand, we note that 
$$
A_{p\over r_0}=\bigcup_{p_0>r_0} A_{p\over p_0}.
$$

\noindent
This implies for all $p>r_0 $ and all $ 
w\in A_{p\over r_0} 
$, the operator $F(L)$ is bounded on $L^p(X, w)$. 
  \hfill{} $\Box$

  \bigskip

\noindent
{\bf Proof of  Theorem~\ref{th3.3}.  }  
 Note that the condition  $\mu(X)<\infty$ implies that $X$ is
  bounded. Hence $X=B(x_0, r_0)$ for some $x_0\in X$ and $0<r_0<\infty$ (\cite{M}).
  It follows from  condition 
   (\ref{e2.3})   that  for any $x\in X,$
$ 
  V(x_0, 1)\leq C  
  \big(1+{d(x, x_0) }\big)^{D}
  V(x, 1 ) 
   \leq 
  C
  V(x, 1 ).
$
This shows that for any $x,y\in X,$
$ 
  \big| K_{e^{-L}}(x,y)\big|  \leq       {C V(x_0, 1)^{-1}}.  
$
 As a consequence,  
\begin{eqnarray}\label{e4.20}
\max\Big\{\big\|e^{-L}\big\|_{L^1(X)\rightarrow L^2(X)}, \, 
  \big\| e^{-L}\big\|_{L^2(X)\rightarrow L^{\infty}(X)}\Big\} \leq C.
\end{eqnarray}

\noindent
On the other hand, for any bounded Borel function $F$ 
such that supp $F\subseteq [0,16]$, the operator
$  F(\sqrt[m]{L}) e^{2L} 
 $ is bounded on $L^2(X)$. This, together with (\ref{e4.20}), yields

\begin{eqnarray*} 
\big\|F(\sqrt[m]{L})\big\|_{L^1(X )\rightarrow   L^{\infty}(X )}
&&
=\big\|e^{-L} \big(F(\sqrt[m]{L}) e^{2L}\big) e^{-L}
\big\|_{L^1(X)\rightarrow   L^{\infty}(X)}\\[3pt]
&&\leq \big\|e^{-L} 
\big\|_{L^1(X)\rightarrow   L^2(X)} 
 \big\|  F(\sqrt[m]{L}) e^{2L}
\big\|_{L^2(X )\rightarrow   L^2(X )}  \big\|e^{-L}
\big\|_{L^2(X)\rightarrow   L^{\infty}(X)}
\\ [3pt]
&&\leq C \|F\|_{L^{\infty}}<\infty.
 \end{eqnarray*}
 
 \noindent
This  implies that  the kernel $K_{F(\sqrt[m]{L})}(x,y)$ of the operator $F(\sqrt[m]{L})$ satisfies
$$
 \sup_{y\in X} \Big|K_{F(\sqrt[m]{L})}(x,y)\Big| \leq C<\infty.
 $$
Hence,   for any $x\in X$,
 \begin{eqnarray*} 
 \big|F(\sqrt[m]{L})f(x)\big|&=& \Big|\int_X K_{F(\sqrt[m]{L})}(x,y) f(y) d\mu(y)\Big|\\
 &\leq& C  \int_X\big|   f(y) \big| d\mu(y)   \\
 &\leq& C\Ma (f)(x),
 \end{eqnarray*}
 
 \noindent
 and  for any $1<p<\infty$ and $w\in A_p$,
 
 $$
 \big\|F(\sqrt[m]{L})f\big\|_{L^p(X, w)} \leq C \big\|\Ma (f)\big\|_{L^p(X, w)} 
 \leq C\big\|f\big\|_{L^p(X, w)}.
  $$

  \bigskip

Therefore, in order to prove Theorem~\ref{th3.3},   
we can assume that supp $F\subset [8, \infty]$. Following  the proof of Theorem 3.1, we set 
 $
F^{\ell}(\lambda)= \varphi(2^{-\ell}\lambda)F(\lambda),
 $
and 
  
  $$
  {\tilde F}=\sum_{\ell \,=3}^{\infty} F^{\ell}\ast \xi,
  $$
  
  \noindent
  where $\xi$ is a function defined in (b) of Lemma~\ref{le4.1}.
  
  By repeating the proof of Theorem~\ref{th3.1} and using (\ref{e4.3}) in place of 
  (\ref{e4.2}) we can prove that
   the operator ${\tilde F}(\sqrt[m]{L})$ is bounded on $L^p(X, w)$ for all $p$ 
   and $w$ satisfying  (i) and (ii) in Theorem~\ref{th3.3}.
 To prove Theorem~\ref{th3.3},  it follows by   Theorem~\ref{th2.4} again that 
 it suffices to show that
for all balls $ B\ni x$,      
\begin{eqnarray}
 \Big(  \oint_{B}
\big|\big(F(\sqrt[m]{L})-{\tilde F}(\sqrt[m]{L})\big)(I-e^{-r_B^mL})^M f \big|^{p_0}d\mu\Big)^{1/p_0} 
\leq C  \Ma \big(  \big| f\big|^{p_0}\big)^{1\over p_0}(x)
\label{e4.21}
\end{eqnarray}
  for all $f\in L^{\infty}_c(X)$.

Let us prove (\ref{e4.21}).  For every   $\ell\geq 3 $ and $r>0$, we set
$H^{\ell}_{r, M}(\lambda)= \big(F^{\ell}(\lambda)- F^{\ell}\ast \xi (\lambda)\big)  
(1-e^{-(r\lambda)^m})^{M},  \lambda>0.$ (Note that  supp $H^{\ell}_{r_B, M}\subseteq [0, 2^{\ell} +1]$.) Now for a given ball $B\subset X$, we put $
 f=\sum\limits_{j=0}^{\infty} f_j,$ where $ f_j=f\chi_{U_j (B)},
 $
and $U_j(B)$ were defined in (\ref{e4.1}).
  We may write
\begin{eqnarray}
\label{e4.22} \hspace{1cm}
\big(F(\sqrt[m]{L})-{\tilde F}(\sqrt[m]{L})\big)  (I-e^{-r_B^mL})^M f  
 &=& 
 \sum_{j=1}^2 \big(F(\sqrt[m]{L})-{\tilde F}(\sqrt[m]{L})\big)  (I-e^{-r_B^mL})^M f_j\nonumber \\
 &+&   
 \lim_{N\rightarrow \infty}\sum_{\ell=3}^{N}  \sum_{j= 3}^{\infty}  
H^{\ell}_{r_B, M}(\sqrt[m]{L})f_j, 
\end{eqnarray}

\noindent
The similar argument as in the proof of Theorem 3.1 gives the desired estimates for $j=1,2.$ 
Next,   fix $j\geq 3$.
For every   $\ell \geq 3$, 
let $K_{H^{\ell}_{r_B, M}(\sqrt[m]{L})}(y,z)$ be  the Schwartz  kernel of operator 
 $H^{\ell}_{r_B, M}(\sqrt[m]{L})$.   Let  ${1\over p_0}-{1\over p_1}={1\over 2},$ 
and denote by   ${1\over p_0}={\theta\over 1} +{1-\theta\over 2}$ and ${1\over p_1}={\theta\over 2},$   
 that is $\theta=2({1\over p_0}-{1\over 2})$. Following  (\ref{e4.10}) and (\ref{e4.11}), we use
  H\"older's inequality and  interpolation again to obtain 
  that for all balls $B\ni x, $     
 \begin{eqnarray}\label{e4.23}
 &&\hspace{-2cm} \Big(  \oint_{B}
\big|H^{\ell}_{r_B, M}(\sqrt[m]{L})  f_j \big|^{p_0}d\mu\Big)^{1/p_0} \nonumber\\ [2pt] 
  &\leq& C2^{jn\over p_0 }   V(B)^{{1\over 2} }  
  \Ma \big(  \big| f\big|^{p_0}\big)^{1\over p_0}(x)  \big\|H^{\ell}_{r_B, M}(\sqrt[m]{L})  \big\|^{1-\theta}_{L^{2}(U_j(B))\to  L^{\infty}(B)}
    \times\\ [2pt] 
	&&\hspace{2cm} \times 
	\big\|\overline{H}^{\ell}_{r_B, M}(\sqrt[m]{L})  \big\|^{\theta}_{L^{2}(B)\to  L^{\infty}(U_j(B))}
 \nonumber 
\end{eqnarray}
 
\noindent
The H\"older  inequality, together with  condition  that $X$ is bounded    give 
\begin{eqnarray}
\label{e4.24} \hspace{1cm}
  \hspace{-2cm}&&\hspace{-1.5cm} \big\|H^{\ell}_{r_B, M}(\sqrt[m]{L})  \big\|^2_{L^{2}(U_j(B))\to  L^{\infty}(B)} \nonumber\\ [3pt]
&=&\sup_{y\in B}  
\int_{U_j(B)}  \Big| K_{H^{\ell}_{r_B, M}(\sqrt[m]{L})}(y,z) \Big|^2
   d\mu(z) \nonumber\\ [3pt]
   &\leq& C  
 \big(2^{j} r_B\big)^{-2s} \sup_{y\in B} \int_{X} \big|K_{H^{\ell}_{r_B, M}(\sqrt[m]{L})}(y,z) 
  \big|^2   d(y, z) ^{2s}d\mu(z)  \nonumber
   \nonumber\\[2pt] 
   &\leq& C_X  
 \big(2^{j} r_B\big)^{-2s} \sup_{y\in B} \int_{X} \big|K_{H^{\ell}_{r_B, M}(\sqrt[m]{L})}(y,z) 
  \big|^2  d\mu(z)   \nonumber
   \nonumber\\[2pt]
  &\leq& \sup_{y\in B}{C_X \over
 V(y, 2^{-\ell})} 
 \big(2^{j} r_B\big)^{-2s}    
 \big\|\delta_{2^{\ell}} \big(H^{\ell}_{r_B, M}\big)\big\|^{2}_{2^{\ell }, q},  
\end{eqnarray}

\noindent
where the last inequality follows from the fact that  supp $H^{\ell}_{r_B, M}\subseteq [0, 2^{\ell} +1]$, 
and then  from (\ref{e3.4})  with 
$N=2^{\ell}$, we have   that
$$
 \int_{X} 
  \big|K_{H^{\ell}_{r_B, M}(\sqrt[m]{L})} (y,z)\big|^2
 d\mu(z) 
  \leq  {C \over
 V(y, 2^{-\ell})}\,
 \big\|\delta_{2^{\ell}} \big(H^{\ell}_{r_B, M}\big)\big\|^{2}_{2^{\ell }, q}.   
$$
From  the expression $H^{\ell}_{r_B, M}(\lambda)=\big(F^{\ell}(\lambda)- F^{\ell}\ast \xi (\lambda)\big) 
  (1-e^{-(r_B\lambda)^m})^{M},$  
one obtains

\begin{eqnarray}
 \|\delta_{2^{\ell}}\big(H^{\ell}_{r_B, M}\big)\|_{2^{\ell }, q}
&=&\big\|\delta_{2^{\ell}}[ F^{\ell}(\lambda)- F^{\ell}\ast \xi (\lambda)\,] (1-e^{-(2^{\ell}r_B\lambda )^m})^{M}
\big\|_{2^{\ell }, q}  \nonumber\\ [4pt]
  &\leq&    
 C {\rm min}\, \big\{1, (2^{\ell}r_B)^{mM}\big\}
  \big\|\delta_{2^{\ell}}[   F^{\ell}(\lambda)- F^{\ell}\ast \xi (\lambda)\,]\big\|_{2^{\ell }, q}.
 \label{e4.25}
\end{eqnarray}
 
 \noindent
Everything then boils down to estimating 
$\|\cdot\|_{2^{\ell}, q}$ norm of 
$\delta_{2^{\ell}}[  F^{\ell}(\lambda)- F^{\ell}\ast \xi (\lambda)\,].$ We make the following 
claim. For its proof, we refer to p. 26, claim (3.29) of \cite{CS} or p. 459, Proposition 4.6 of  \cite{DOS}.
 
 \bigskip
 
 \begin{prop}\label{prop4.2}
 Suppose that $\xi\in C_c^{\infty}$ is a function such that 
 {\rm supp} $\xi\subset [-1, 1]$, $\xi\geq 0$, ${\hat \xi}(0)=1$,
 ${\hat \xi}^{(\kappa)}(0)=0$ for all $1\leq \kappa\leq [s]+2$ and set $\xi_N(t)=N\xi(Nt)$. 
  Assume also that
 {\rm supp} $G\subset [0,1]$. Then
 
 $$
 \big\|G-G\ast \xi_N \big\|_{N, q}\leq C N^{-s}\big\|G\big\|_{W^q_s}
 $$
 
 \noindent
 for all $s>1/q.$
 \end{prop}

 \bigskip
 
 \noindent
By Proposition~\ref{prop4.2} 
\begin{eqnarray*} 
\big\|\delta_{2^{\ell}}[ F^{\ell}(\lambda)- F^{\ell}\ast \xi (\lambda)\,]\big\|_{2^{\ell},\, q}=
\big\|\delta_{2^{\ell}}[ \varphi_{\ell} F] 
- \xi_{2^{\ell}}\ast \delta_{2^{\ell}}[ \varphi_{\ell} F]\big\|_{2^{\ell}, q}
\leq C 2^{-\ell s}
\big\|\delta_{2^{\ell}}[\varphi_{\ell} F]\big\|_{W^q_s},
\end{eqnarray*}
 \noindent
 and thus
 
 \begin{eqnarray}
 \|\delta_{2^{\ell}}\big(H^{\ell}_{r_B, M}\big)\|_{2^{\ell}, q}
  &\leq&    
 C 2^{-\ell s} {\rm min}\, \big\{1, (2^{\ell}r_B)^{mM}\big\} 
 \big\|\delta_{2^{\ell}}[\varphi_{\ell} F]\big\|_{W^q_s}.
 \label{e4.26}
\end{eqnarray}

\noindent
Substituting (\ref{e4.26})  back into  (\ref{e4.24}), 
we then use the doubling property (\ref{e2.2})
to obtain

\begin{eqnarray*}
&&\hspace{-1cm}  \big\|H^{\ell}_{r_B, M}(\sqrt[m]{L})  \big\|^2_{L^{2}(U_j(B))\to  L^{\infty}(B)} \\[3pt]
&\leq& {C\over V(B)}  2^{-2s j} 
   \big(2^{\ell} r_B\big)^{-2s}   {\rm min}\, \big\{1, (2^{\ell}r_B)^{2mM}\big\}
 \max\big\{1, \big(2^{\ell} r_B\big)^n\big\}   
\big\|\delta_{2^{\ell}}[\varphi_{\ell} F]\big\|_{W^q_s}^2,
\end{eqnarray*}
 
\noindent
which is exactly the same estimate as in (\ref{e4.16}).
 
 \medskip
 
Following the proof of Theorem 3.1, an argument as above shows the same  
 estimate (\ref{e4.17}) for the term  $\big\|\overline{H}^{\ell}_{r_B, M}(\sqrt[m]{L})  \big\|_{L^{2}(B)\to  L^{\infty}(U_j(B))}.$
The rest of the proof of   (\ref{e4.21})  is just a
  repetition of   the proof of  Theorem~\ref{th3.1}, 
    so we skip it.
  Hence, we complete the proof of Theorem~\ref{th3.3} when $X$
has  a finite measure, i.e., $\mu(X)<\infty$.
  \hfill{} $\Box$

\bigskip

 \section{Two interpolation results}
\setcounter{equation}{0}

In this section we continue to assume that 
$L$ is a non-negative self-adjoint  operator on $L^2({X})$, which 
has a kernel $p_t(x,y)$ satisfying a Gaussian upper bound $(GE)$.
Using interpolation, other conditions on the weight can be found which guarantee that $F(L)$
is a bounded operator. We first prove the following result.

\medskip

\begin{theorem}\label{th5.1}  
Let  $s>\frac{n}{2} $ and let 
$r_0=\max\{{2(n+{D})\over 2s+{D}},1\}$. 
Suppose that    the operator $L$ satisfies condition {\rm (\ref{e3.2})}   
  with some   $q\in [2, \infty].$ 
If $1<p<\infty$ and $w^{r_0}\in A_p$, 
 then for any bounded Borel function $F$ such that $\sup_{t>0}\|\eta\, \delta_tF\|_{W^q_s}<\infty,$
 the operator $F(L)$ is bounded on $L^p(X, w)$. Moreover,
\begin{eqnarray*} 
   \|F(L)  \|_{L^p(X, w)\rightarrow L^p(X, w)}\leq    C_s\Big(\sup_{t>0}\|\eta\, \delta_tF\|_{W^q_s}
   + |F(0)|\Big).
\end{eqnarray*}
\end{theorem}

  \noindent
{\bf Proof.}
We will derive  Theorem~\ref{th5.1} from  Theorem~\ref{th3.1} by using Proposition~\ref{prop2.4}
 and the characterization of $A_p$ functions
that if $w\in A_p$, then there are $A_1$ weights $u$ and $v$ such that $w=u v^{1-p}$ 
(\cite{J}).

Following the proof of Theorem 2 of \cite {KW}, we fix $p, 1<p<\infty,$ and $w$ so that $w^{r_0}\in A_p$ where $r_0={2(n+{D})\over 2s+{D}}$. 
We have that $w^{r_0} =u v^{1-p}$, $u, v\in A_1$,
or $w=u^{r_0^{-1}} v^{1-p\over r_0}$.  Next, write this as

$$
w=u^{r_0^{-1}} v^{1-p\over r_0} =\big(u^{\alpha} v^{\beta}\big)^{t} 
\big(u^{\gamma} v^{\delta}\big)^{1-t}=w_0^t  w_1^{1-t}.
$$

\noindent
in which 

\begin{eqnarray}
\label{e5.1}
\alpha t +\gamma (1-t)= r_0^{-1},
\end{eqnarray}
\begin{eqnarray}
\label{e5.2}
\beta t +\delta (1-t)= r_0^{-1} (1-p).
\end{eqnarray}

\noindent
Then in order to use Proposition~\ref{prop2.4} for weights which satisfy Theorem~\ref{th3.1},
we require

\begin{eqnarray}
\label{e5.3}
w_0^{-{1\over r-1}}\in A_{r^{'}\over r_0}, \ \ \ \ 1<r<\min\big\{ r_0^{'}, p \big\}
\end{eqnarray}
\begin{eqnarray}
\label{e5.4}
w_1 \in A_{q \over r_0}, \ \ \ \ q>\max\big\{ r_0, p \big\}.
\end{eqnarray}
\begin{eqnarray}
\label{e5.5}
t={q-p\over q-r}.
\end{eqnarray}

\noindent
Recall that $u\in A_1$ (similarly $v\in A_1$) implies

$$
 \oint_B u \leq C u(x)\ \ \ \ \mbox{for almost all $x\in B$}.
$$

\noindent
Therefore, if $\alpha>0$ and $\beta<0$, letting $s={r^{'}\over r_0}$, we have

\begin{eqnarray*}
 &&\hspace{-2cm}\Big(  \oint_B  w_0 ^{-{1\over r-1}}  \Big)  
\Big(  \oint_B  w_0^ {1\over (r-1)(s-1)}  \Big)^{s-1}\\ [3pt]
 &\leq&  \Big(  \oint_B  u^{-{\alpha\over r-1}}\, v^{-{\beta\over r-1}}  \Big)  
\Big(  \oint_B  u^{ {\alpha\over (r-1)(s-1)}} \, v^{{\beta\over (r-1)(s-1)}} \Big)^{s-1}  \\ [3pt]
&\leq& \Big(  \oint_B  u \Big)^{-{\alpha\over r-1}}
  \Big(  \oint_B  v^{-{\beta\over r-1}}  \Big)   
\Big(  \oint_B  v \Big)^{{\beta\over r-1}} 
 \Big(  \int_B  u^{ {\alpha\over (r-1)(s-1)}}   \Big)^{s-1}  \\ [3pt]
 &\leq& C,
\end{eqnarray*}

\noindent
if 
$$
\alpha=(r-1)\Big({r^{'}\over r_0}-1\Big)= {r\over r_0}-r+1\ \ \ {\rm and}\ \ \beta=-(r-1);
$$

\noindent
this is $w_0^{-{1\over r-1}}\in A_{r^{'}\over r_0}$ for these
 values of $\alpha$ and $\beta$. Similarly, we can show
$w_1\in A_{q\over r_0}$ if $\gamma=1$ and $\delta=-\big({q\over r_0}-1\big)$.
 Using these values of $\alpha$ and $\gamma$,
we have (\ref{e5.1}) if $t={1\over r.}$ Next, solving (\ref{e5.2}) for $q$, we get $q=r'(p-1)$. This value of $q$
also satisfies (\ref{e5.5}). Therefore, if we choose $r< \min\big\{ r_0^{'}, p\big\}$
 close enough to $1$ so  that $q=r'(p-1)>\max\big\{r_0, p\big\}$, then (\ref{e5.1})- (\ref{e5.5}) hold.
 This proves
Theorem~\ref{th5.1}.
 \hfill{} $\Box$

\bigskip

If $X={\Bbb R}^n$ then Theorem~\ref{th5.1} can be strengthen for the following polynomial weights. When $w(x)
=|x|^{\beta}$, we have $w\in A_p$ if $-n<\beta<n(p-1)$. 
Applying  Theorem~\ref{th3.1} and 
Theorem 3 of  \cite{KW} to such $w$ and using interpolation with change of measures, we have
the following theorem.
 
\medskip

\begin{theorem}\label{th5.2} Let  $s>\frac{n}{2} .$ 
Suppose that the operator $L$ satisfies condition (\ref{e3.2})  
  with some   $q\in [2, \infty].$ 
If $1<p<\infty$ and $\max\{-n, -sp\}<\beta<\min\{n(p-1), sp\}$, 
 then for any bounded Borel function $F$ such that $\sup_{t>0}\|\eta\, \delta_tF\|_{W^q_s}<\infty,$
 the operator $F(L)$ is bounded on $L^p({\Bbb R}^n, |x|^{\beta})$. In addition,
\begin{eqnarray*} 
   \|F(L)  \|_{L^p({\Bbb R}^n, |x|^{\beta})\rightarrow L^p({\Bbb R}^n, |x|^{\beta})}\leq  
     C_s\Big(\sup_{t>0}\|\eta\, \delta_tF\|_{W^q_s}
   + |F(0)|\Big).
\end{eqnarray*}

\noindent
In particular, if   $s<n$  and ${ n \over  s}<p<\big({ n \over  s}\big)^{'},$ we get $-n<\beta<n(p-1)$; 
we may also take $p={ n \over  s}$ and $p=\big({ n \over  s}\big)^{'}$.
\end{theorem}

\medskip

\noindent {\bf Proof.} The proof of Theorem~\ref{th5.2} can be obtained by making minor modifications
with the proof of Theorem 3 in \cite{KW} and using the Theorem~\ref{th3.1}.  We give a brief argument
of this proof for completeness and convenience for the reader.

 Notice that $-n\geq -sp$\,  if\,  $n/s\leq p$, and $n(p-1)\leq sp$ \, if \, $p\leq (n/s)'$. 
 Therefore, for $s<n$ the conclusion of Theorem 
 5.2 can be divided into three cases:
 
 \begin{eqnarray}
\label{e5.6}
1<p<{n\over s}\ \ {\rm and}\ \ -sp<\beta<n(p-1),
\end{eqnarray}
\begin{eqnarray}
\label{e5.7}
 {n\over s}\leq p\leq  \Big({n\over s}\Big)'  \ \ {\rm and}\ \ -n<\beta <n(p-1),
\end{eqnarray}
\begin{eqnarray}
\label{e5.8}
  \Big({n\over s}\Big)' <p<\infty \ \ {\rm and}\ \ -n<\beta <sp.
\end{eqnarray}
 
 \noindent
 Since (\ref{e5.8}) is the dual of (\ref{e5.6}), we need only concern ourselves with (\ref{e5.6}) and (\ref{e5.7}).

 Because $|x|^{\beta}\in A_p$ if and only if $-n<\beta<n(p-1)$, it follows from Theorem~\ref{th3.1}  that for $s<n,$ $F(L)$ is bounded 
 on $L^p({\Bbb R}^n, |x|^{\beta})$ if

 \begin{eqnarray}
\label{e5.9}
 {n\over s} \leq p<\infty \ \ {\rm and}\ \ -n<\beta < ps-n,
\end{eqnarray}
\begin{eqnarray}
\label{e5.10}
  1<p\leq \Big({n\over s}\Big)'   \ \ {\rm and}\ \ -n+p(n-s)<\beta <n(p-1).
\end{eqnarray}
 
 \noindent
 However, combining (\ref{e5.9}) and (\ref{e5.10}), we have (\ref{e5.7}) and are left with only proving 
 (\ref{e5.6}).
 
 Let $q=n/s$ and $r<n/s;$ then also $r<(n/s)'$. By (\ref{e5.10}) and (\ref{e5.7}), $F(L)$ is bounded on 
 $L^p({\Bbb R}^n, |x|^{\beta_0})$ and $L^p({\Bbb R}^n, |x|^{\beta_1})$ for $-n+r(n-s)<\beta_0<n(r-1)$ 
 and $-n<\beta_1<n(q-1)$. Using Proposition~\ref{prop2.4}, if $r<p<q$ we see that $F(L)$ is bounded on 
$L^p({\Bbb R}^n, |x|^{\beta})$ for 

$$
\beta=\beta_0\Big({q-p\over q-r}\Big)  +\beta_1\Big({p-r\over q-r}\Big).
$$

\noindent
Thus $\beta$ satisfies

$$
\big\{-n +r(n-s)\big\}\Big({q-p\over q-r}\Big) -n \Big({p-r\over q-r}\Big) <\beta <n(r-1)\Big({q-p\over q-r}\Big) +n(q-1)\Big({p-r\over q-r}\Big).
 $$
 Simplifying and using the fact that $q=n/s$, we get
 
 \begin{eqnarray}
\label{e5.11}
 {n^2(r-1)\over n-sr} +{psr(s-n)\over n-sr} <\beta<n(p-1).
\end{eqnarray}
 But, as $r\rightarrow 1$, the left-hand side of (\ref{e5.11})  tends to $-sp.$ So, taking $r$
 sufficiently close to $1$ allows us to choose any $\beta$ satisfying $-sp<\beta<n(p-1)$.
 
 When $s=n$, the restriction in Theorem~\ref{th5.2} is $-n<\beta<n(p-1)$ for $1<p<\infty.$ But, when $s=n$
 in Theorem~\ref{th3.1}, we require $w\in A_p$, and $|x|^{\beta}\in A_p$ if $-n<\beta<n(p-1)$.
 This completes the proof of Theorem~\ref{th5.2}. 
 \hfill{} $\Box$

\bigskip

Note that in the case of Fourier spectral multipliers  
  Theorem~\ref{th5.2} is best possible, except for endpoint equalities 
for $\beta $  see \cite[pp. 360-361]{KW}.

\bigskip

 \section{Applications}
\setcounter{equation}{0}

\subsection{Homogeneous groups}\label{sec61}

Let ${\bf G}$  be a Lie group of polynomial growth and let $X_1, ..., X_k$ be a system of 
 left-invariant vector fields on ${\bf G}$ satisfying the H\"ormander condition. We define
 the Laplace operator $L$ acting on $L^2({\bf G})$ by the formula
 \begin{eqnarray}
L=-\sum_{i=1}^k X_i^2.
\label{e6.1}
\end{eqnarray}  

\noindent If $B(x,r)$ is the ball define by the distance associated with system $X_1, ..., X_k$
(see e.g. Chapter III.4, \cite{VSC}), then there exist natural numbers $n_0, n_{\infty}\geq 0$
such that $V(x, r) \sim r^{n_0}$ for $r\leq 1$ and  $V(x, r)\sim r^{n_{\infty}}$
for $r>1$ (see e.g. Chapter III.2, \cite{VSC}). Note that this implies that doubling condition \eqref{e2.2} holds with 
the doubling dimension  $n=\max\{n_0,{n_{\infty}}  \}$. Note also that one can take $D=0$ 
in the estimates \eqref{e2.3}.   We call ${\bf G}$  a homogeneous group
if there exists a family of dilations on ${\bf G}$. A family of dilations on a Lie group
${\bf G}$ is a one-parameter group $({\tilde\delta }_t)_{t>0}$  $({\tilde\delta }_t \circ {\tilde\delta }_t 
={\tilde\delta }_{t s})$  of automorphisms of ${\bf G}$ determined by
 \begin{eqnarray}
{\tilde\delta }_t Y_j =t^{n_j} Y_j ,
\label{e6.2}
\end{eqnarray} 

\noindent
where $Y_1, ..., Y_{\ell}$ is a linear basis of Lie algebra of ${\bf G}$ and $n_j\geq 1$
for $1\leq j\leq \ell$ (see \cite{FS}). We say that an operator $L$  defined by (\ref{e6.1})
 is homogeneous if  ${\tilde\delta }_t X_i =t X_i $ for $1\leq i\leq k$ and the system 
 $X_1, ..., X_k$ satisfies the H\"ormander condition. Then for the sub-Riemannian geometry 
 corresponding to the system   $X_1, ..., X_k$  one has $n_0=n_{\infty}=\sum_{j=1}^{\ell} n_j$ (see \cite{FS}).
 Hence the doubling dimension is equal to $n=n_0=n_{\infty}$. 
 
 Spectral multiplier theorems for the homogeneous Laplace operators acting on homogeneous groups were
 investigated  by Hulanicki and Stein \cite{HS}, Folland and Stein  \cite[Theorem 6.25 ]{FS},  and De Michele 
 and Mauceri \cite{DeM}.  See also \cite{C} and \cite{MM}.  We have the following weighted spectral multiplier result. 

\begin{prop}\label{prop6.1} Let $L$ be the homogeneous sub-Laplacian 
defined by the formula {\rm (\ref{e6.1})} acting on
a homogeneous group ${\bf G}$.  Then Theorem~{\rm \ref{th3.1}} holds for spectral multipliers $F(L)$  with 
$q=2$, $D=0$ and the doubling dimension given by $n=n_0=n_{\infty}$.
\end{prop}
 \noindent {\bf Proof.}  It is well known that the heat kernel  corresponding 
to the operator $L$ satisfies Gaussian bound $(GE)$. 
It is also not difficult to check that for some constant $C>0$
$$
\|F(\sqrt L)\|_{2 \to \infty}^2 =C \int_0^\infty |F(t)|^2 t^{n-1}dt.
$$
See for example equation (7.1) of \cite{DOS} or \cite[Proposition 10]{C}.
It follows from the above equality that  the operator $L$ satisfies estimate (\ref{e3.2})
with $q=2$. Hence Theorem~\ref{th3.1} holds for spectral multipliers $F(L)$ with $q=2$, $D=0$ and $n=n_0=n_{\infty}$.
\hfill{}$\Box$

   \bigskip

This result can be extended to ``quasi-homogeneous'' operators  acting on homogeneous groups, see \cite{S1} and \cite{DOS}.

 \medskip

In the setting of general Lie groups of polynomial growth,  spectral multipliers were
investigated by Alexopoulos. Our following weighted spectral multiplier result extends
 Alexopoulos's  unweighted result in \cite{A1}.

\begin{prop}\label{prop6.2} Let $L$ be a group invariant  operator
 acting on
a Lie group ${\bf G}$ of polynomial growth defined by (\ref{e6.1}). Then  Theorem ~\ref{th3.2} 
holds for spectral multipliers $F(L)$ with the doubling dimension $n=\max\{n_0,{n_{\infty}}  \}$ and $D=0$.
 \end{prop}

 \noindent {\bf Proof.}  It is well known that the heat kernel  corresponding 
to the operator $L$ satisfies Gaussian bound $(GE)$ so  the operator $L$ satisfies  estimate (\ref{e3.2})
for $q=\infty$, see the proof of Theorem~\ref{th3.2} above and Lemma 2.2 of \cite{DOS}. Hence Theorem ~\ref{th3.2} holds for spectral multipliers $F(L)$.
   \hfill{}$\Box$

   \bigskip

\bigskip

\subsection{Compact manifolds}

For a general non-negative self-adjoint elliptic operator on a compact manifold,
 Gaussian bound $(GE)$  holds  by general elliptic regularity theory. Further, one has
the Avakumovi${\rm{\breve c}}$-Agmon-H\"ormander theorem.

\medskip

\begin{theorem} \label{th6.3}
Let $L$ be a non-negative elliptic pseudo-differential operator
of order $m$ on a compact manifold $X$ of dimension $n$. Then
\begin{eqnarray}\label{e6.3}
\big\|\chi_{[R, R+1]} (L^{1/m})\big\|^2_{L^1(X)\rightarrow L^2(X)}\leq C R^{n-1}, \ \ \ 
\forall \, R\in{\Bbb R}^+.
\end{eqnarray}
\end{theorem}

\medskip

Theorem~\ref{th6.3} was proved by H\"ormander \cite{Ho2}.
This theorem has the following useful consequence.

\medskip

\begin{cor} \label{cor6.4}
Condition {\rm (\ref{e3.4})} with $q=2$   holds for non-negative
  elliptic pseudo-differential operators on 
  compact manifolds.
\end{cor}

\medskip

\noindent
{\bf Proof:} By spectral theorem
\begin{eqnarray*} 
\sup_{y\in X} \big\|K_{F(\sqrt[m]{L})}(\cdot, y)\big\|^2_{L^2(X) }&\leq& \Big( 
\sum_{\ell=1}^N\big\|\chi_{[\ell-1, \, \ell]} F(\sqrt[m]{L})\big\|^2_{L^1(X)\rightarrow L^2(X)}\Big)^{1/2}\\
&\leq& C N^{n/2} \big\|\delta_N F\big\|_{N,2}
\end{eqnarray*}

\noindent
as required.
  \hfill{} $\Box$

  \medskip
  
The importance of estimate (\ref{e6.3}) for multiplier theorems was noted by Sogge \cite{Sog1}, 
who used it to establish  the convergence of the Riesz means up to the critical exponent
$(n-1)/2$ (see also \cite{CSo}).

\medskip

 \begin{prop}\label{prop6.5}  Suppose that 
  $L$ is a  non-negative self-adjoint   elliptic differential  operator
  of order $m\geq 2$ acting  on a compact Riemannian manifold $X$ of dimension $n.$  
  Then the operator $L$ satisfies estimate {\rm (\ref{e3.4})}
for $q=2$, and hence Theorem~\ref{th3.3}  holds for spectral multipliers $F(L)$ under the same 
conditions with $q=2$, $D=0$ and with the doubling dimension $n$. That is the exponent $n$ in 
{\rm \eqref{e2.2}} is equal to the topological dimension of the
manifold $X$. 
 \end{prop}

 \medskip
  
  \noindent {\bf Proof.}  This result is a direct consequence of Theorem~\ref{th3.3} and
 Corollary~\ref{cor6.4}.
   \hfill{}$\Box$
 
\bigskip

Proposition~\ref{prop6.5} applied to an elliptic operator on a compact Lie group gives a stronger
result than  Proposition~\ref{prop6.2}. One can say that for elliptic operators on a compact
Lie group  Proposition~\ref{prop6.1} holds. However, we do not know if the
 Avakumovi${\rm{\breve c}}$-Agmon-H\"ormander   condition holds for sub-elliptic operators on a compact
 Lie group (see also \cite{CS}). Hence, Proposition~\ref{prop6.2} gives the strongest known
 result for sub-elliptic operators on a compact Lie group.
 
 \bigskip
 
 \subsection{Laplace operators on irregular domains with Dirichlet boundary conditions}
 
 Let $\Omega$ be a connected open subset of ${\Bbb R}^n.$ Note that if the boundary of
 $\Omega$ is not smooth enough, then
 $\Omega$ is not necessarily a homogeneous space because the doubling condition
 might not hold. 
 
 In this section we are interested in dealing with weighted norm estimates in those contexts.
 As it is pointed out in \cite{DM}, one can extend the singular operators defined in $\Omega$
 to the space ${\mathbb R}^n.$ Since there is no assumption on the regularity of the kernels in 
 space variables, the extension of the kernel still satisfies similar conditions. Given $T$, a bounded linear operator
 on $L^p(\Omega), 1<p<\infty$, the extension of $T$ to ${\mathbb R}^n$ is defined as 
 ${\widetilde T}f(x)=T\big(f\chi_{\Omega}\big)(x)\chi_{\Omega}(x)$ for $f\in L^p({\Bbb R}^n)$. Then, $T$ is bounded
 on $L^p(\Omega)$ if and only if ${\widetilde T}$ is bounded  on $L^p({\Bbb R}^n).$
 If $K$ is the kernel of $T$, then the associated kernel of ${\widetilde T}$  is given by
  ${\widetilde K}(x,y)= K(x,y)$ for $(x,y) \in \Omega\times \Omega$ and  ${\widetilde K}(x,y)= 0$ otherwise. As it is observed in \cite{DM},
  the assumptions on the kernels do not involve their regularity so they imply similar 
  properties on the kernels of the extended operators.
 
 We are going to use the notation $A_p({\Bbb R}^n)$ in order to make clear that the Muckehhoupt 
 weights are considered in the whole space ${\Bbb R}^n.$ 
 The following result gives examples of singular
 integral multipliers on spaces without the doubling condition.

  \begin{prop}\label{prop6.6} Suppose that  $\Delta_\Omega$ is the Laplace operator with Dirichlet
  boundary condition $\Omega\subset {\Bbb R}^n$. 
 Let $s>{n/2}$ and    $r_0=  \max\big(1, {n/s}\big)$.
 Then for any bounded Borel function $F$
 such that 
$\sup_{t>0}\|\eta\, \delta_tF\|_{W^{\infty}_s}<\infty,$
 the operator $F(\Delta_\Omega)$ is bounded on $L^p(\Omega, w)$ for all $p$ and $w$ satisfying
  $  r_0  < p<\infty$ and   $w\in A_{p/r_0}({\Bbb R}^n)$. \noindent  
In addition,
\begin{eqnarray*} 
   \| F(\Delta_\Omega) \|_{L^p(\Omega, w)\rightarrow L^p(\Omega, w)}\leq    C_s\Big(\sup_{t>0}\|\eta\, \delta_tF\|_{W^{\infty}_s}
   + |F(0)|\Big).
\end{eqnarray*}
 
 \end{prop}

 \medskip
  
  \noindent {\bf Proof.}  Note that 
  $$
  0\leq  K_{\exp(-t\Delta_{\Omega})}(x,y)\leq {1\over (4\pi t)^{n/2}} \exp\Big(-{|x-y|^2\over 4t}\Big)
  $$
  
  \noindent
  (see e.g., Example 2.18, \cite{Da1}).
  That is  the heat kernels  corresponding to $\Delta_{\Omega}$
 satisfy Gaussian bounds  $(GE)$, and the operator $\Delta_\Omega$ satisfies estimate {\rm (\ref{e3.2})}
for $q=\infty$. Then, Proposition~\ref{prop6.6} follows from estimate (\ref{eff}) 
and Theorem~\ref{th3.2}  applied to the extended operator ${\widetilde {F(\Delta_\Omega)}}$.   
Hence the same weighted norm estimates hold for the original operator
${ {F(\Delta_\Omega)}}$.
 \hfill{}$\Box$

   \bigskip
   
   \subsection{Schr\"odinger operators}
   
   In this section we discuss applications of our main results to spectral multipliers of Schr\"odinger operators.

  Let  $\Delta$ be the standard Laplace operator acting on ${\Bbb R}^n $. We consider the Schr\"odinger operator 
   $L=-\Delta+V$ where $V: {\Bbb R}^n\rightarrow {\Bbb R}, V\in L^1_{\rm loc}({\Bbb R}^n)$ and
   $V\geq 0.$ The operator $L$ is defined by the quadratic form. If
   $p_t(x,y)$ denotes the heat kernel corresponding to $L$ then as a consequence of the Trotter
   product formula
 \begin{eqnarray}\label{e6.4}
   0 \le  p_t(x,y)\leq \tilde p_t(x,y),
  \end{eqnarray}

   \noindent
   where $\tilde p_t(x,y)$ denotes the standard Gauss  heat kernel corresponding to $\Delta$. 
   
   The estimate  (\ref{e6.4}) holds also for heat kernel $p_t(x,y)$ of 
   Schr\"odinger operator with electromagnetic potentials, see \cite[Theorem 2.3]{Si}
   and \cite[(7.9)]{DOS}. For the Schr\"odinger operator in this setting, 
 estimate (\ref{e3.2}) holds for $q=\infty$ as in the next result.
   
   \medskip

 \begin{prop}\label{prop6.7}  Assume that  
   $L=-\Delta+V$ where $\Delta$ is the standard Laplace operator acting on ${\Bbb R}^n $  and  
   $  V\in L^1_{\rm loc}({\Bbb R}^n )$ is a non-negative function.
    Then the operator $L$ satisfies estimate {\rm  (\ref{e3.2})}
for $q=\infty$, and hence  Theorem~{\rm \ref{th3.2}} holds for spectral multipliers $F(L)$ under the same 
conditions with $q=\infty$, $D=0$ and the doubling constant $n$.  
 \end{prop}

We note that under suitable  additional assumptions this result
  can be extended by a similar proof to situation of magnetic Schr\"odinger operators
  acting on a complete Riemannian manifold with non-negative potentials.

 \medskip
  
  \noindent {\bf Proof.}  This result is a consequence of (\ref{e6.4}) and  Theorem~\ref{th3.2}.
   \hfill{}$\Box$
 
 \medskip

\subsection{Estimates on operator norms of holomorphic functional calculi}

For $\theta>0$, we put $\sum_{\theta}=\{z\in {\bf C}-\{0\}: |{\rm arg} \, z|<\theta \}$.
Let $F$ be a bounded holomorphic function on $\sum_{\theta}.$ By $\|F\|_{\theta, \infty}$
we denote the supremum of $F$ on $\sum_{\theta}.$ We are interesting in finding sharp bounds,
in terms of $\theta$, of the norm of $F(L)$ as the operator acting on $L^p(X, w)$. 
The following proposition, which is a weighted version of   \cite[Proposition 8.1]{DOS},
 is a consequence of Theorem~\ref{th3.2}.
 
\medskip

\begin{prop}\label{prop6.8}  Let $L$ be an operator satisfying assumptions of Theorem {\rm \ref{th3.2}}
Let $s>{n\over 2}$ and let   $r_0=  \max\big\{1, {2(n+{D})\over 2s+{D}}\big\}$.
 Then 
 the operator $F(L)$ is bounded on $L^p(X, w)$ for all $p$ and $w$ satisfying
  $  r_0  < p<\infty$ and   $w\in A_{p\over r_0}$. 
  In addition,
\begin{eqnarray*} 
   \|F(L)  \|_{L^p(X, w)\rightarrow L^p(X, w)} 
     \leq {C_{\epsilon}\over \theta^{{n \over 2}+\epsilon} } \big\|F\big\|_{\theta, \infty} 
\end{eqnarray*}

\noindent
for every $\epsilon>0,$   $  r_0  < p<\infty$ and   $w\in A_{p\over r_0}.$
\end{prop}  

 \medskip
  
  \noindent {\bf Proof.}  
It is easy to check, using the Cauchy formula that there exists a constant $C$
independent of $F$ and $\theta$ such that

\begin{eqnarray*}
\sup_{\lambda>0}\big|\lambda^k F^{(k)}(\lambda)\big|\leq {C\over \theta^k}
 \big\|F\big\|_{\theta, \infty},  \ \ \ \ \ \forall k\in {\Bbb Z}_+.
 \end{eqnarray*}
 
 \noindent
 For any $\epsilon>0$, $\sup_{t>0}\big\|\eta\delta_t F\big\|_{W^{\infty}_{k-\epsilon}} 
 \leq C\sup_{\lambda>0}\big|\lambda^k F^{(k)}(\lambda)\big| $  so by  interpolation
 
\begin{eqnarray*}
\sup_{t>0}\big\|\eta\delta_t F\big\|_{W^{\infty}_{s}} 
 \leq {C_{\epsilon}\over \theta^{s+\epsilon}} \big\|F\big\|_{\theta, \infty}. 
  \end{eqnarray*}
  
  \noindent
  Applying the above inequality and Theorem~\ref{th3.2} we obtain 
  Proposition~\ref{prop6.8} (see also, Theorem 4.10, \cite{CDMY}). 
   \hfill{}$\Box$
 
 \medskip

\vskip 1cm

 \noindent
{\bf Acknowledgments.}  The authors are grateful to the referee who gave detailed
comments and suggestions for improving the original manuscript.
This work was started during   the third named  author's  
stay at  Macquarie University.  
  L.X. Yan would like to thank   the Department of Mathematics of Macquarie University
for its hospitality.  The research of X.T. Duong was  supported by  
Australia Research Council (ARC).  
 The research of L.X. Yan was supported by   Australia Research Council (ARC)   
and  NNSF of China   (Grant No.  10771221 and 10925106).

\end{document}